\documentclass[letterpaper, 10 pt, conference]{ieeeconf}  

\IEEEoverridecommandlockouts                              

\def\Nset{\mathbb{N}}
\def\Rset{\mathbb{R}}
\def\cR{{\cal R}}
\newcommand{\comment}[1]{}

\def\diag{{\mathrm {\tt diag}}}

\usepackage{graphics} 
\usepackage{times} 

\usepackage{mathtools,amsmath,amssymb}
\usepackage[pdftex]{hyperref}
\usepackage{color}

\newtheorem{theorem}{Theorem}

\newtheorem{remark}{Remark}
\newtheorem{definition}{Definition}

\title{\LARGE \bf
Feedback Control Principles for Biological Control of Dengue Vectors*}

\author{Pierre-Alexandre Bliman$^{\star}$
\thanks{* A shorter version of this report has been accepted for presentation at the {\em European Control Conference} ECC19, Naples (Italy), June 25-28, 2019.
\newline
\indent
The author acknowledges the financial support of the Franco-Columbian program ECOS-Nord, through the project C17M01.
This work benefited from fruitful discussions with Y. Dumont (CIRAD, University of Pretoria) within the framework of the Phase 2 SIT feasibility project against \textit{Aedes albopictus} in Reunion Island, jointly funded by the Regional Council of La Reunion and the European Regional Development Fund (ERDF) under the 2014-2020 Operational Program.}
\thanks{$^{\star}$ Pierre-Alexandre Bliman is with Sorbonne Universit\'e, Universit\'e Paris-Diderot SPC, Inria, CNRS, Laboratoire Jacques-Louis Lions, \'equipe Mamba, F-75005 Paris, France; e-mail: \href{mailto:pierre-alexandre.bliman@inria.fr}{pierre-alexandre.bliman@inria.fr}}%
}

\begin{document}

\maketitle
\thispagestyle{empty}
\pagestyle{empty}

\begin{abstract}

Controlling diseases such as dengue fever, chikungunya and zika fever by introduction of the intracellular parasitic bacterium Wolbachia in mosquito populations which are their vectors, is presently quite a promising tool to reduce their spread.
While description of the conditions of such experiments has received ample attention from biologists, entomologists and applied mathematicians, the issue of effective scheduling of the releases remains an interesting problem for Control theory.
Having in mind the important uncertainties present in the dynamics of the two populations in interaction, we attempt here to identify general ideas for building release strategies, which should apply to several models and situations.
These principles are exemplified by two interval observer-based feedback control laws whose stabilizing properties are demonstrated when applied to a model retrieved from \cite{Bliman:2018aa}.
Crucial use is made of the theory of monotone dynamical systems.
\end{abstract}

\section{Introduction}
\label{se1}

Already a major issue in southern countries since decades, the control of vector-borne diseases transmitted by mosquitoes recently became an important concern for northern countries too, due to the geographical dissemination of the vectors favored by climate change, urbanization and increasing international travel.
When no vaccine or curative treatment exists, traditional control measures focus on reducing the vector population in order to reduce the pathogen transmission.
Mechanical control methods by preventive actions and removal of breeding sites are invaluable, but costly and difficult to implement efficiently.
Chemical control by insecticides has been quite popular, but on top of its negative impacts on humans and animals, it experiences decreasing efficiency due to resistance generation.
Recently, various biological control methods have been proposed and tested as alternative or complementary strategies, typically by the release of transgenic or sterile mosquitoes \cite{Alphey:2010aa,Alphey:2014aa}.
A new promising strategy based on the use of  the bacterium {\em Wolbachia}, is being developed to control the spread of  dengue fever and other diseases transmitted by species of the genus {\em Aedes} (chikungunya, zika fever, yellow fever).
It is grounded in the fact that {\em Wolbachia} severely reduces the insect vectorial ability,  indirectly by reducing lifespan and fertility, and directly by reducing the virus ability to proliferate within the organism \cite{McMeniman:2009aa,Moreira:2009aa,Hoffmann:2011aa,Wilder-Smith:2017aa}.

The dynamics resulting from the introduction of {\em Wolbachia}-infected mosquitoes in wild populations has been abundantly considered, even in the absence of spatial aspects; see \cite{Keeling:2003aa,Farkas:2010aa,Zheng:2014aa,Zheng:2017aa,Yakob:2017aa,Xue:2017aa,Campo-Duarte:2017aa,Campo-Duarte:2018aa,Bliman:2018aa,Almeida:2018aa}, and \cite{Hughes:2013aa,Koiller:2014aa,Ndii:2015aa} for the effects induced on a dengue epidemics.
Field sampling methods allow to evaluate population density \cite{Focks:2004aa,Silver:2007aa}, and such operations are commonly performed during experimental campaigns.
Also, the presence of the bacterium in the captured sample may be investigated by polymerase chain reaction (PCR) method \cite{Hoffmann:2011aa}.
%
Theoretically, this opens up the possibility to assess the released quantities on the basis of the observed population, and to benefit from the multiple advantages of closed-loop methods.
Most papers analyze only the effect of a unique, initial, release.
By contrast, \cite{Zheng:2017aa} considers continuous releases and \cite{Yakob:2017aa} impulsive releases (with no analytical results), both with constant release amplitude, while \cite{Bliman:2018aa} provides linear control-law and \cite{Campo-Duarte:2017aa,Campo-Duarte:2018aa,Almeida:2019aa} optimal control approach.
However these stabilization results are highly dependent upon the setting, and their generalization is in no way straightforward. 
Our aim here is to identify simple control principles, based on the understanding of the biological phenomena involved,
and to test them on the model introduced in \cite{Bliman:2018aa}.
More specifically, we propose two control laws achieving complete infestation, respectively through the introduction of infected adults and larvae.
Their implementation necessitates the construction of an interval observer \cite{Moisan:2009aa,Efimov:2013aa,Efimov:2016aa} for the considered system.
%

The content of the paper is as follows.
The characteristics of the compartmental model developed in \cite{Bliman:2018aa} are recalled in Section \ref{se2}, together with some analysis results.
Section \ref{se3} proposes a class of interval observers for this model.
The main results (Theorems \ref{th4} and \ref{th5}) are exposed and proved in Section \ref{se4}.
Numerical simulations are shown in Section \ref{se6}, and some concluding remarks are given in Section \ref{se7}.

{\bf Notations.}
The following notations are used for the positive, resp.\ negative, parts of a real number $z$: $|z|_+ := \max\{z,0\}$, $|z|_- := -\min\{z,0\}$.
Both are nonnegative, and
\begin{equation}
\label{eq26}
\forall z\in\Rset,\qquad z=|z|_+-|z|_-\ .
\end{equation}

The usual order relation $\geq$ in $\Rset$ is employed, and the same notation is extended to vectors and matrices by the product order: $(x_1, x_2)\geq 0$ {\em iff} $x_1\geq 0$ and $x_2\geq 0$.
It is similarly extended to locally integrable functions taking on values in one of these ordered spaces, with the evident meaning.
In all cases, the symbol $>$ is used as usual to mean `at most equal to, but not equal'.
As an example, for locally integrable functions $f,g$ defined on a common domain, `$f>g$' is equivalent to: `$f\geq g$ and $f>g$ on a nonzero-measured set'.
Also, we use the notation $\ll$ for `much less than'.
Whenever the sign $\pm$ is used, it means that the relevant formula is valid for $+$ and for $-$.
Last, usual matrix notations are employed: $I_n, 0_n$ denote respectively the identity and null square matrices of order $n$, $n\in\Nset$, $^{\mbox{\tiny\sf T}}$ denotes transposition.

\section{Modeling}
\label{se2}

\subsection{Model equations}
\label{se21}

The compartmental model in \cite{Bliman:2018aa} focuses on the main effects pertinent for control purposes.
It contains two life phases: a preliminary phase which gathers the early, aquatic, stages (egg, larva and pupa), subject to competition for food and space; and an adult, aerial, mature phase.
The corresponding state variables are respectively denoted $L$ (`larvae') and $A$ (`adults'), and the uninfected and {\em Wolbachia}-infected populations are distinguished by indices $U$ and $W$.
Assuming in first approximation that the sex ratio is stationary and the mortality is sex-independent, these variables
represent indifferently the quantities of males or females, up to constant ratio.
After normalization, the $4$-dimensional population model used in the present paper is the following \cite{Bliman:2018aa}:
\vspace{-.1cm}
\begin{subequations}
\label{eq1}
\begin{gather}
\label{eq1a}
\dot L_U = \gamma_U{\cal R}_0^U \frac{ A_U}{ A_U+ A_W} A_U - (1+ L_U+ L_W) L_U \\
\label{eq1b}
\dot A_U =  L_U -\gamma_U  A_U\\
\label{eq1c}
\dot L_W = \gamma_W {\cal R}_0^W  A_W - (1+ L_U+ L_W) L_W +u_L(t) \\
\label{eq1d}
\dot A_W = L_W -\gamma_W  A_W + u_A(t)
\end{gather}
\end{subequations}

The quantities $\gamma_U, \gamma_W$ are (normalized) mortality rates, while $\cR_0^U, \cR_0^W$ are the {\em basic offspring numbers} \cite{Yang:2009aa,Ferreira:2014aa} associated to each population.
They represent the average number of mosquitoes born to each adult mosquito during its entire lifespan.

The signal $u_L$ (resp.\ $u_A$) represents the infected larvae (resp.\ adults) released per time unit in order to infect the wild population.
By construction both have nonnegative values.
Using $x := (L_U,A_U,L_W,A_W)$, $u:=(u_L,u_A)$ as state and input variables, the previous controlled system writes compactly:
\begin{equation}
\label{eq00}
\dot x = f(x) + Bu,\qquad y=Cx\ ,
\end{equation}
where the function $f$ is easily expressed from \eqref{eq1}, and the input matrix $B$ is given by
\[
B:= \begin{pmatrix}
0 & 0 & 1 & 0\\
0 & 0 & 0 & 1
\end{pmatrix}^{\mbox{\tiny\sf T}}\ .
\]
The output matrix $C$ is used to define which are the available measurements, it is specified in Section \ref{se3}.

Well-posedness of the Cauchy problem for system \eqref{eq1} does not present specific difficulty, e.g.\ for locally integrable input $u$ with nonnegative values almost everywhere.
As a result of this and of the fact that $(f(x))_i\geq 0$ whenever $x\geq 0$ and $x_i=0$, $i=1,\dots,4$, the set $\Rset_+^4$ is forward invariant for the flow.

\subsection{Phenomenological hypotheses on the infection dynamics}
\label{se22}

The main phenomenological assumptions contained in the previous infection model are the following.

\noindent $\circ$ 1.\
The {\em Wolbachia} infection in {\em Aedes} mosquito leads to lifespan shortening through adult mortality increase \cite{McMeniman:2009aa,Yeap:2011aa} and fecundity rate reduction \cite{Ruang-areerate:2006aa}.
This is accounted for by the following inequalities between normalized constants:
\begin{equation}
\label{eq12}
0<\gamma_U < \gamma_W,\qquad \cR_0^U>\cR_0^W>1\ .
\end{equation}
The variations of larval development time and mortality are regarded as limited and are not modeled here.

\noindent $\circ$ 2.\
The density-dependent mortality is a major component of larval dynamics \cite{Southwood:1972aa,Focks:1993aa}, rendered here by quadratic competition terms in the immature phase dynamics.
Accounting for restricted food and space in the breeding sites, it acts evenly on the immature infected and uninfected population, with a rate proportional to its global size $L_U+L_W$.

\noindent $\circ$ 3.\
The {\em Wolbachia} transmission occurs only vertically, from mother to offspring.
It is accompanied by a phenomenon called {\em cytoplasmic incompatibility}, which provides a reproduction advantage to infected females against uninfected ones and facilitates the spread of the bacterium.
Cytoplasmic incompatibility is characterized by the fact that an uninfected female only produces offspring when mating with an uninfected male \cite{Influential:1998aa,McMeniman:2009aa,Yeap:2011aa,Walker:2011aa}.
On the contrary when an infected female mates, the outcome is infected independently of whether the male is infected or not.
The probability of a male being uninfected is $\frac{A_U}{A_U+A_W}$, and the nonlinear birth term $\frac{A_U}{A_U+A_W}A_U$ in \eqref{eq1a} accounts for such encounters.
By contrast, the birth term for infected in \eqref{eq1c} is simply proportional to the corresponding population $A_W$.


\subsection{Analysis of the uncontrolled model \rm\cite{Bliman:2018aa}}
\label{se23}

In order to understand the meaning of the basic offspring numbers, consider the following auxiliary system:
\begin{equation}
\label{eq13}
\dot L = \gamma\cR A-(1+L)L,\qquad \dot A = L-\gamma A\ ,
\end{equation}
with $\gamma>0, \cR>0$.
System \eqref{eq13} is monotone \cite{Smith:1995aa}.
When $\cR<1$, $0_2$ is the unique (extinction) equilibrium and is globally asymptotically stable; while when $\cR>1$ it is unstable, and the nonzero equilibrium $(\cR-1,\frac{\cR-1}{\gamma})$ appears, which is asymptotically stable and limit point of every nonzero trajectory.
The relevance of this comes from the following fact, which may be easily checked: when initially one of the two populations (infected or uninfected) is absent from \eqref{eq1}, then the other one evolves alone according to system \eqref{eq13}.

Transposing this remark discloses immediately the existence, under hypothesis \eqref{eq12}, of two equilibria for system \eqref{eq1} with zero input $u$, characterized by the state values $(L_U^*,A_U^*,0,0)$ and $(0,0,L_W^*,A_W^*)$, where:
\begin{equation*}
L_\eta^* := \cR_0^\eta-1,\quad A_\eta^* := \frac{L_\eta^*}{\gamma_\eta},\qquad \eta=U,W\ .
\end{equation*}
They correspond respectively to {\em disease-free} and {\em complete infestation} situations.
Under hypothesis \eqref{eq12}, both are locally asymptotically stable: this assumption ensures the sustainability of each of the two isolated populations, with a greater basic offspring number for the uninfected population, in accordance with the
fact that {\em Wolbachia} globally reduces the fitness of the infected mosquitoes.

Two other equilibria exist, which correspond to {\em extinction} and {\em coexistence}, and are both unstable.
System \eqref{eq1} presents {\em bistability}, and the control purpose is typically to pass from the disease-free situation to the complete infestation.

A structural property, central to analyze the behavior of system \eqref{eq1}, is {\em monotonicity} \cite{Smith:1995aa}.
It will be instrumental in the sequel.
We first introduce a specific order relation defined in the space $\Rset^4$.
\begin{definition}
\em
The ordering induced in the space $\Rset^4$ by the cone $\Rset_-\times\Rset_-\times\Rset_+\times\Rset_+$ is denoted $\succeq$: by definition,
\begin{equation}
\label{eq27}
x\succeq x'\ \Leftrightarrow\ x_i\leq x'_i,\ i=1,2 \text{ and } x_i\geq x'_i,\ i=3,4
\end{equation}
\end{definition}

\noindent
The following results is borrowed from \cite{Bliman:2018aa}.
\begin{theorem}
\label{th0}
\em
System \eqref{eq1} with zero input $u$ is {\em strongly order-preserving} {\rm\cite[p.\ 2]{Smith:1995aa}} in $\Rset_+^4$ for the order induced by \eqref{eq27}.
\end{theorem}

\section{Upper and lower state estimates}
\label{se3}

In order to study control synthesis in Section \ref{se4}, upper and lower estimates of the state variables will be needed.
Notice that using monotonicity and the fact that the input variable $u$ takes on nonnegative values, one may show e.g.\ that the solutions of \eqref{eq1} fulfill, for any initial condition, the properties:
\begin{gather*}
0\leq \liminf_{t\to +\infty} L_U(t)\leq \limsup_{t\to +\infty} L_U(t) \leq L_U^*,\\
0\leq \liminf_{t\to +\infty} A_U(t)\leq \limsup_{t\to +\infty} A_U(t) \leq A_U^*
\end{gather*}
The use of such rough estimates is possible, but may yield
uselessly costly control laws.
On the other hand, it is important to stress here that the stabilisation results that we present below in Section \ref{se4} do {\em not} require any convergence property from the used observers.

We assume the availability of, say, $p$ measurements, and introduce the output variable
$y=Cx$
for some fixed matrix $C\in\Rset^{p\times 4}$.
Typically $y$ contains measurement on the values of $L_U(t), L_W(t)$ and/or $A_U(t), A_W(t)$, depending upon the information monitored in the field.

In order to take into account measurement uncertainties, one assumes that are available two (locally integrable, nonnegative-valued) auxiliary signals $y_-(t), y_+(t)$ bounding (componentwise) the exact, but unknown, value $y(t)$: for (almost) any $t\geq 0$,
$0 \leq y_-(t) \leq y(t) \leq y_+(t)$.

As a quite easy consequence of the monotonicity of the model (see Theorem \ref{th0}), one has the following result.

\begin{theorem}\em (Interval observer design).
\label{th2}
For $x$ evolving according to \eqref{eq00}, let the evolution of the variables
\vspace{-.1cm}
\begin{equation}
\label{eq28}
x_- := \begin{pmatrix}
L_U^+ \\ A_U^+ \\ L_W^- \\ A_W^-
\end{pmatrix},\qquad
x_+ := \begin{pmatrix}
L_U^- \\ A_U^- \\ L_W^+ \\ A_W^+
\end{pmatrix}
\end{equation}
be defined by
\vspace{-.1cm}
\begin{subequations}
\label{eq18}
\begin{gather}
\label{eq18a}
\hspace{-.25cm}
\dot x_- = f(x_-(t)) + Bu +K_-(x_-(t))(y_-(t)-Cx_-(t))\\
\label{eq18b}
\hspace{-.25cm}
\dot x_+ = f(x_+(t)) + Bu +K_+(x_+(t))(y_+(t)-Cx_+(t))
\end{gather}
\end{subequations}
for gain matrices $K_-(\cdot), K_+(\cdot)\in\Rset^{4\times p}$ depending continuously upon the state.
Assume that
\vspace{-.1cm}
\begin{subequations}
\label{eq17}
\begin{gather}
\label{eq17a}
\diag\{ -I_2; I_2\}  K_\pm(\cdot) \geq 0\\
\label{eq17b}
\diag\{ -I_2; I_2\}  K_\pm(\cdot)C\ \diag\{ -I_2; I_2\}\leq 0
\end{gather}
and that the $i$-th line of $K_-(\cdot)$ (resp.\ $K_+(\cdot)$) is null whenever the $i$-th component of its argument is null, $i=3,4$ (resp.\ $i=1,2$).
\end{subequations}
If
\begin{equation}
\label{eq19}
x_-(t) \preceq x(t) \preceq x_+(t)
\ \text{ and }\
0 \leq x_-(t),x(t),x_+(t)
\end{equation}
for $t=0$, then the same holds true for any $t\geq 0$.
\end{theorem}

Recall that $f$ in the previous statement is the function that allows to write \eqref{eq1} as \eqref{eq00}.
Theorem \ref{th2} provides sufficient conditions under which, if
\[
0\leq L_\eta^-(t) \leq L_\eta(t) \leq L_\eta^+(t),\
0\leq A_\eta^-(t) \leq A_\eta(t) \leq A_\eta^+(t),
\]
$\eta=U,W$, for $t=0$, then the same holds true for any $t\geq 0$.
In addition, it is easy to show that if $x_\pm(t)=x(t)$, say for $t=0$, then the same is true for any $t\geq 0$.
System \eqref{eq18} therefore constitutes an {\em interval observer}\footnote{Notice however that this does not mean that the estimates are bounded, as according to the input $u$, the state $x$ itself may grow unbounded.} \cite{Moisan:2009aa,Efimov:2013aa,Efimov:2016aa} for \eqref{eq1} --- but {\em not} an {\em asymptotic observer}, as previously mentioned.

Coming back to the definition of the order $\preceq$ in \eqref{eq27}, conditions \eqref{eq17a}-\eqref{eq17b} appear as rephrasing of the conditions:
\begin{subequations}
\label{eq16}
\vspace{-.3cm}
\begin{multline}
\label{eq16a}
\forall y,y'\in\Rset^p,\ y \leq y'\ \Rightarrow\\ K_-(\cdot)y\preceq K_-(\cdot)y',\ K_+(\cdot)y\preceq K_+(\cdot)y'
\end{multline}
\vspace{-.9cm}
\begin{multline}
\label{eq16b}
\forall x,x'\in\Rset^4,\ x\preceq x'\ \Rightarrow\\ K_-(\cdot)Cx \succeq K_-(\cdot)Cx',\ K_+(\cdot)Cx \succeq K_+(\cdot)Cx'
\end{multline}
\end{subequations}

\begin{remark}
\label{re0}
Notice that \eqref{eq17b} comes as a consequence of \eqref{eq17a} when e.g.
\vspace{-.1cm}
\addtocounter{equation}{-3}
\begin{subequations}
\addtocounter{equation}{2}
\begin{equation}
\label{eq17c}
C\ \diag\{ -I_2; I_2\}\leq 0\ .
\end{equation}
When for example measurements of the two larval populations $L_U(t), L_W(t)$ are available, then $p=2$ and one may take $C= \begin{pmatrix}
1 & 0 & 0 & 0\\ 0 & 0 & -1 & 0
\end{pmatrix}$, in accordance with \eqref{eq17c}.
Condition \eqref{eq17a} then says that the first two rows of the admissible gain matrices $K_-(\cdot), K_+(\cdot)$ must be nonpositive, and their last two rows nonnegative.
\hfill$\square$
\end{subequations}
\addtocounter{equation}{2}
\end{remark}

\begin{remark}
\label{re1}
A simple way to fulfill the last technical assumption on the matrix gain lines is to introduce a smoothed Heaviside function for the corresponding components; namely to take the $i$-th line of $K_\pm(x_\pm)$ equal to the $i$-th line of a constant matrix $M_\pm$ fulfilling \eqref{eq16}, multiplied by the function $\max\{\min\{x_{\pm,i},\varepsilon\},0\}$, for some $0<\varepsilon \ll 1$.
\hfill$\square$
\end{remark}

{\em Proof of Theorem {\rm \ref{th2}}.}
For simplicity, we omit in the sequel the argument of the matrix-valued functions $K_\pm$.

When $x_-\preceq x$ and $y_-\leq y=Cx$, one has, using \eqref{eq16}:
$K_-(y_--Cx_-) = K_-(y_--y) + K_-C(x-x_-)
\preceq 0$;
and similarly, when $x_+\succeq x$ and $y_+\geq y$:
$K_+(y_+-Cx_+) = K_+(y_+-y) + K_+C(x-x_+)
\succeq 0$.
Therefore, for system \eqref{eq00}-\eqref{eq18} one has
$\dot x_- \preceq f(x_-(t))$, $\dot x_+ \succeq f(x_+(t))$.

Use of Kamke's theorem \cite{Coppel:1965aa} then allows to show that when the first property in \eqref{eq19} holds for $t=0$, it holds also for any $t\geq 0$, provided that all components of the three vectors $x, x_-, x_+$ remain nonnegative --- otherwise the composition by $f$ is not licit, and the solution is not defined.

The nonnegativity of $x$, whose evolution does not depend upon $x_-, x_+$, has already been established \cite{Bliman:2018aa}.
We now show that the values of $x_-, x_+$ effectively remain nonnegative along time.
First, one verifies easily that the indexes mentioned in the statement are exactly those related to the evolution of $L_U^-, A_U^-, L_W^-, A_W^-$.
On the other hand, when $(x_\pm)_i=0$ for some $i=1,\dots,4$, the $i$-th component $(f(x_\pm))_i$ is nonnegative.
Thus, for any matrix line concerned by the condition imposed in the statement, $(\dot x_\pm)_i\geq 0$ whenever $(x_\pm)_i=0$, and for this reason the values of the corresponding signals never leave the interval $[0,+\infty)$.
With an initialization done in accordance with \eqref{eq19}, the lower estimates are thus nonnegative, which forces in turn the components $L_U^+, A_U^+, L_W^+, A_W^+$ to be nonnegative.
Therefore the second inequality in \eqref{eq19} also holds, and this achieves the proof of Theorem \ref{th2}.
\hfill $\blacksquare$

\section{Feedback control principles and observer-based stabilization}
\label{se4}

Due to the uncertainty inherent to the models of population dynamics, it seems valuable to propose feedback control strategies obeying {\em general}, `model-free', principles.
We propose in the sequel two such approaches, related to control by release of adult mosquitoes and of larvae.

\subsection{Control by release of adult mosquitoes}
\label{se41}

First, it is clear that, if the inter-species competition induced by the presence of infected mosquitoes is sufficient to make unviable the uninfected population, then the latter disappears.
Adults do not directly interact, as shown by equations \eqref{eq1b} and \eqref{eq1d}.
However, they participate centrally to the natality.
With this in mind, assume that, due to releases of adults through the input $u_A$ in \eqref{eq1d}, the fractional term responsible in \eqref{eq1b} for the cytoplasmic incompatibility, is kept smaller than $\frac{1}{\cR_0^U}$, where $\cR_0^U$ is the basic offspring number of the uninfected population.
Under such conditions, the uninfected population evolves sensibly as in equation \eqref{eq13} --- however with a {\em subcritical} `apparent' basic offspring number, as $\cR_0^U\times \frac{1}{\cR_0^U}=1$.
This suggests that targeting the cytoplasmic incompatibility term may be a way to realize eradication of the uninfected population.
This general principle underpins the following result.

\begin{figure*}
{\renewcommand{\theequation}{\bf A}
\begin{equation}
\label{eq8b}
u_L \geq 0,\
u_A \geq 
K
\begin{pmatrix}
L_U^+\\
A_U^+\\
A_U^-\\
L_W^-\\
A_W^+\\
A_W^-
\end{pmatrix}
\text{ on } [T,+\infty),\quad
K := \begin{pmatrix}
k_U & k_U|k -\gamma_U|_+ & k_U|k -\gamma_U|_- & -1 &  |\gamma_W -k|_+ &  |\gamma_W -k|_-
\end{pmatrix}
\end{equation}
}
\hrule
\end{figure*}

\begin{figure*}
\centering
\includegraphics[width=6cm]{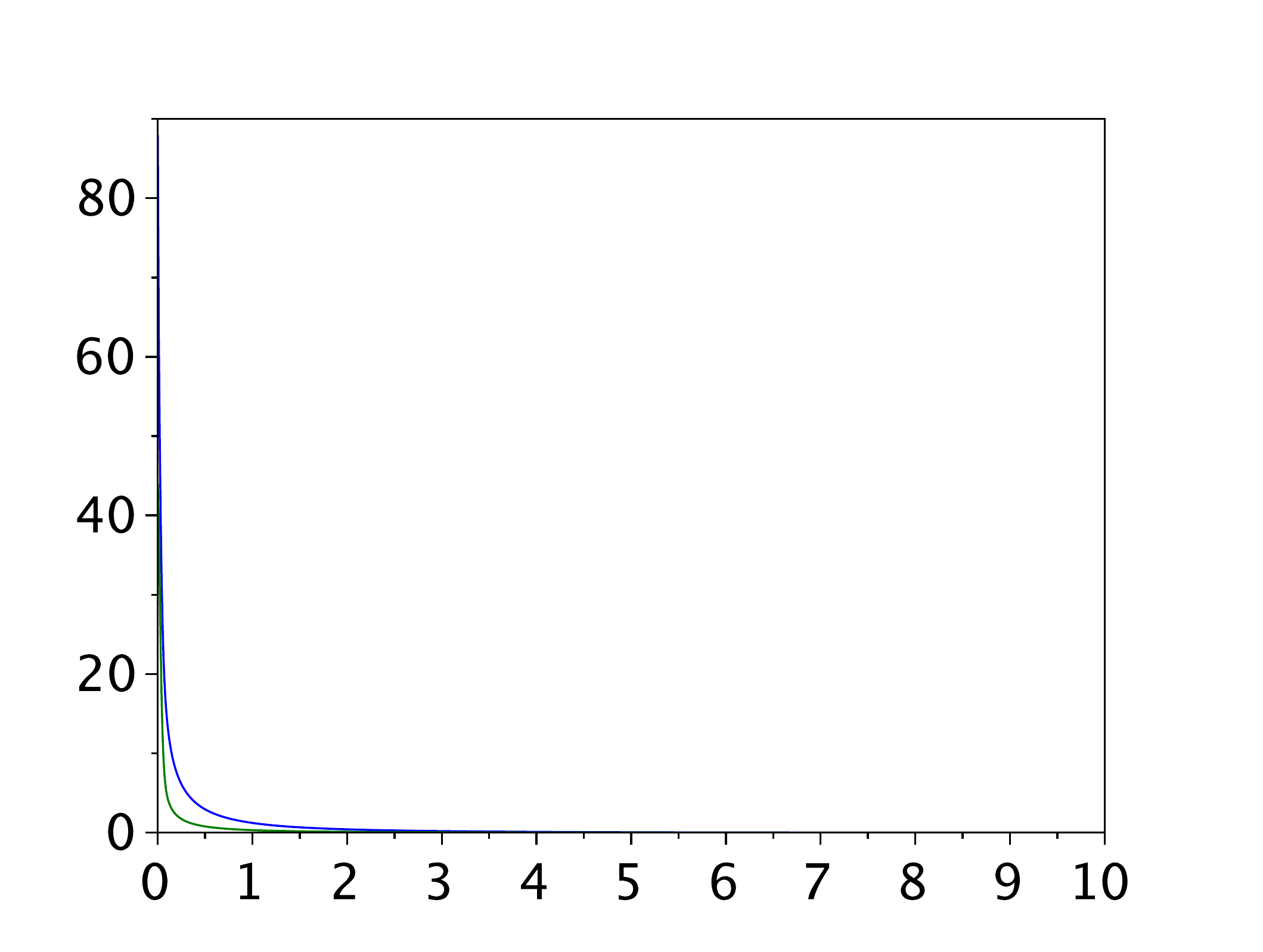}
\includegraphics[width=6cm]{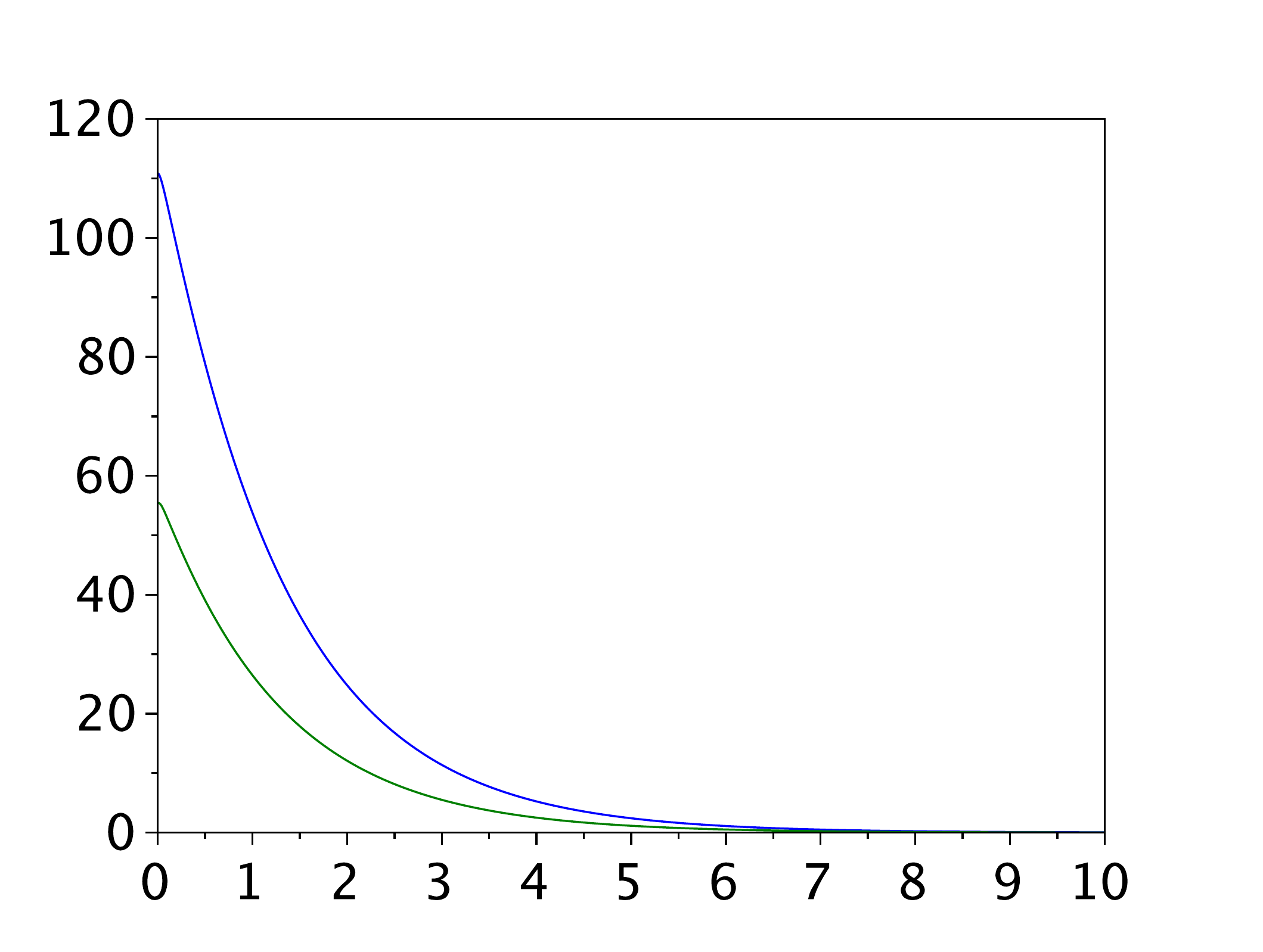}\\
\includegraphics[width=6cm]{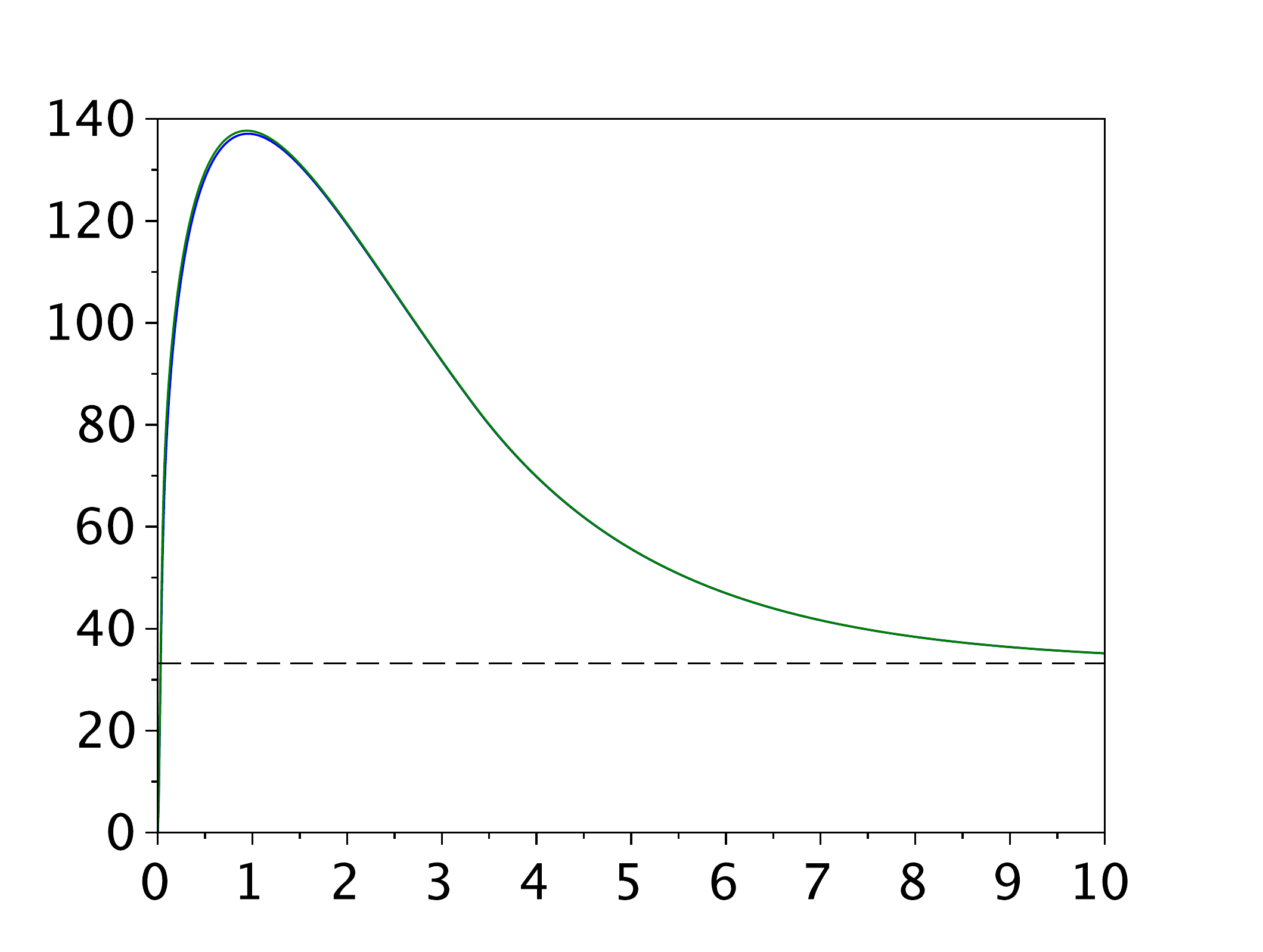}
\includegraphics[width=6cm]{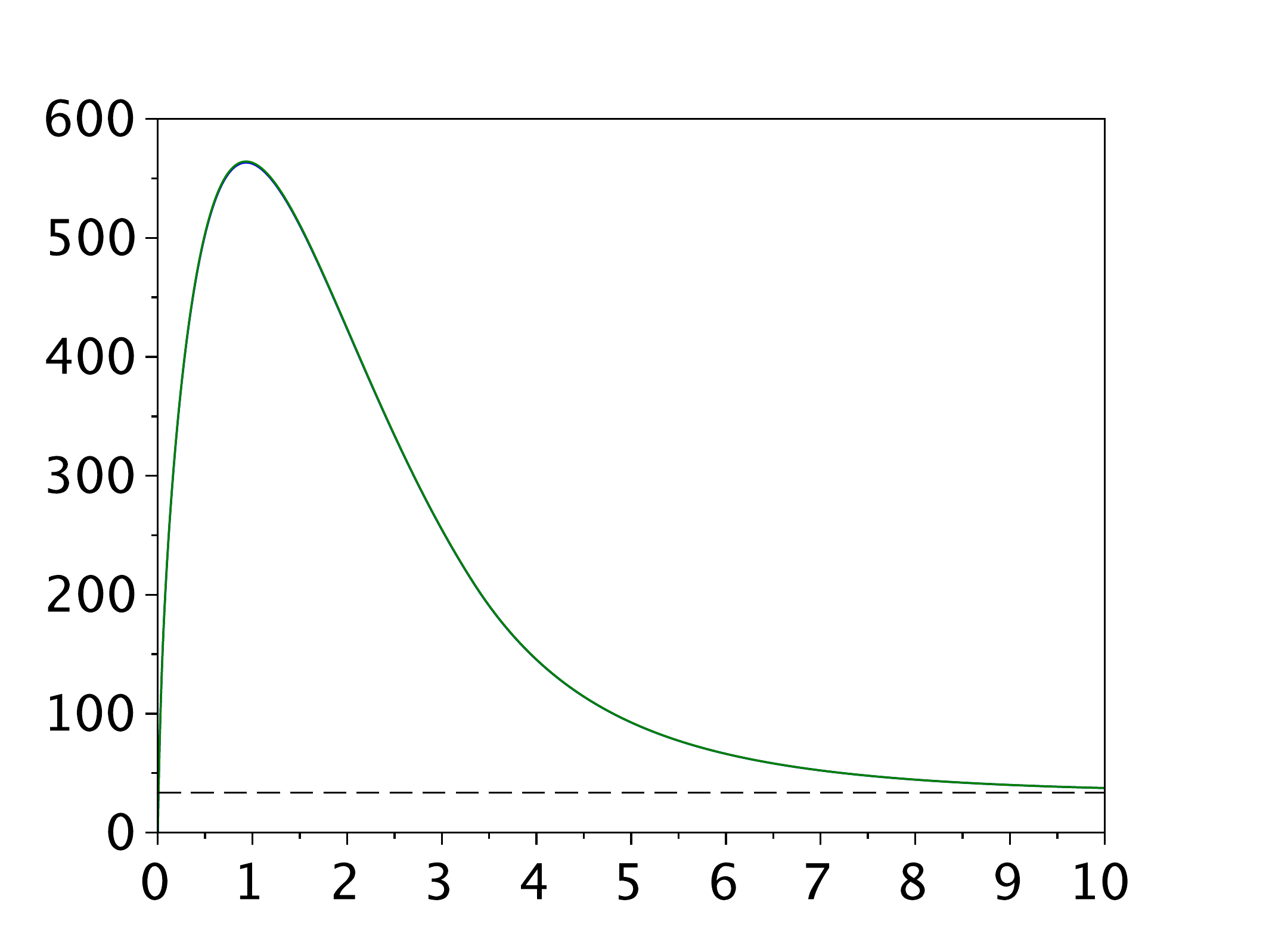}
\caption{{\bf Release of adult mosquitoes}: evolution of the uninfected (top) and {\em Wolbachia}-infected (bottom), as function of time.
The larvae appear on the left column, the adults on the right one.
The components of the state $x$ (resp.\ of the estimate $x_-$) appear in {\color{green}\bf green} (resp.\ in {\color{blue}\bf blue}).}
\label{fig1}
\end{figure*}

\begin{figure*}
\centering
\includegraphics[width=6cm]{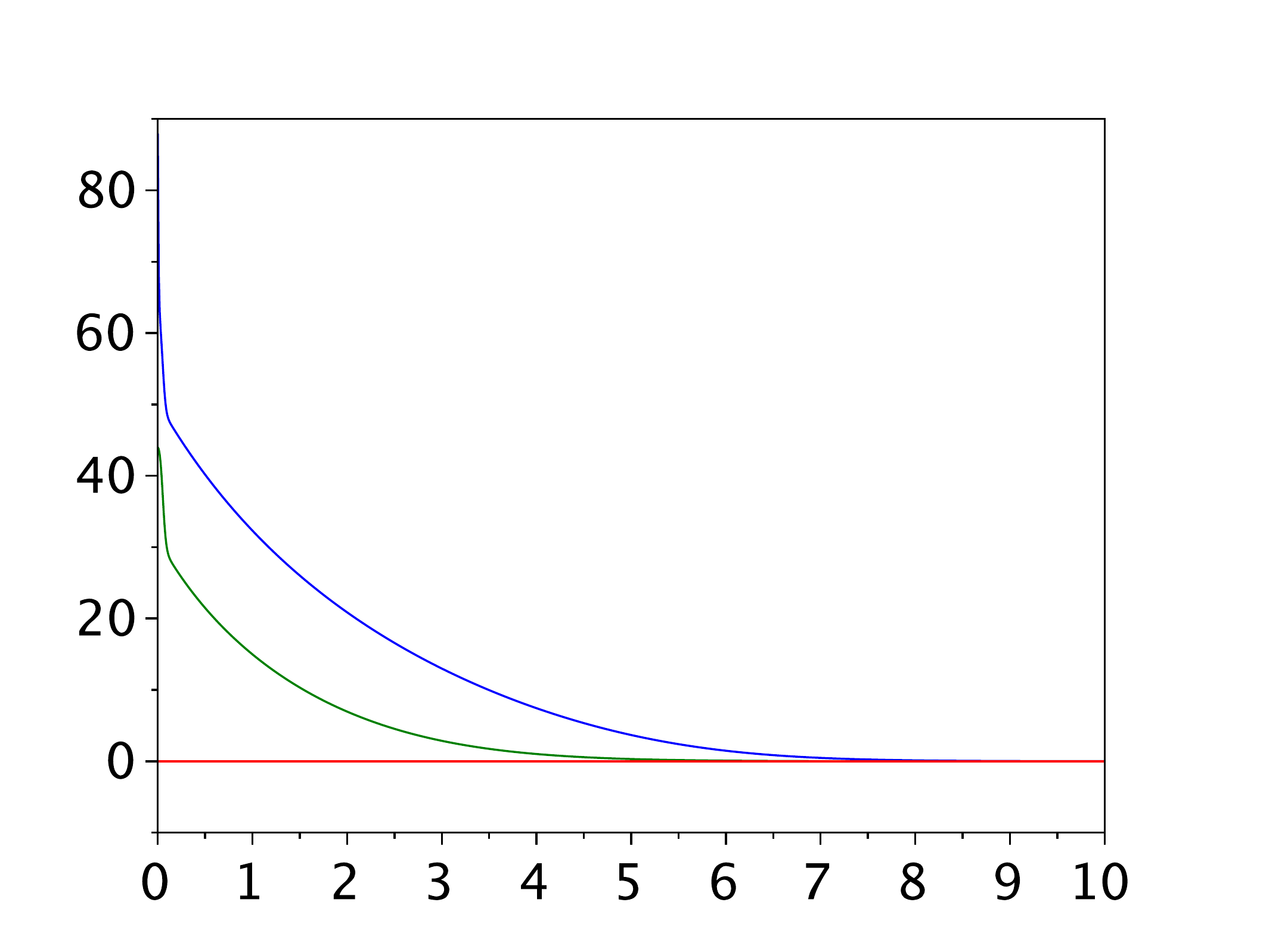}
\includegraphics[width=6cm]{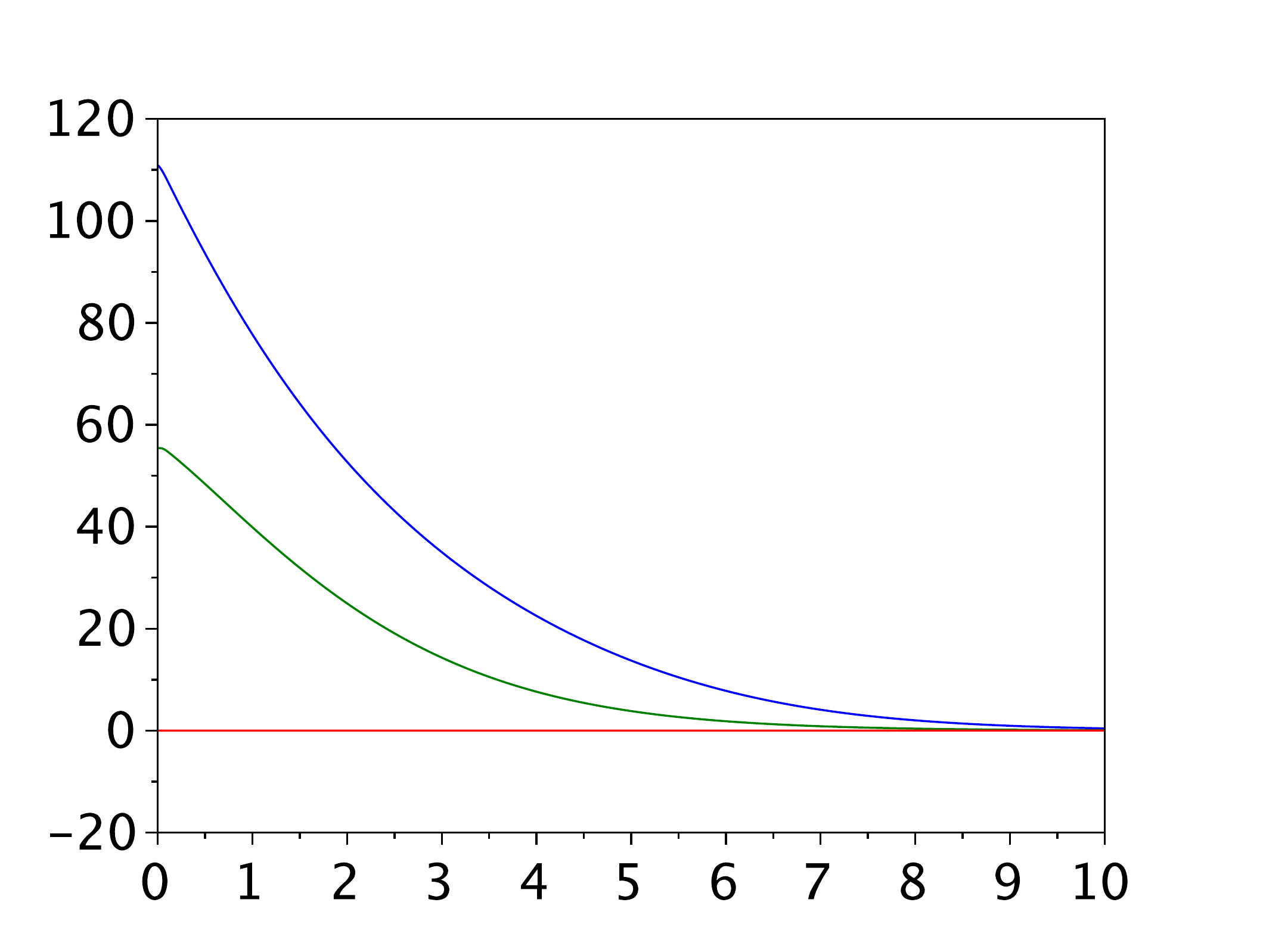}\\
\includegraphics[width=6cm]{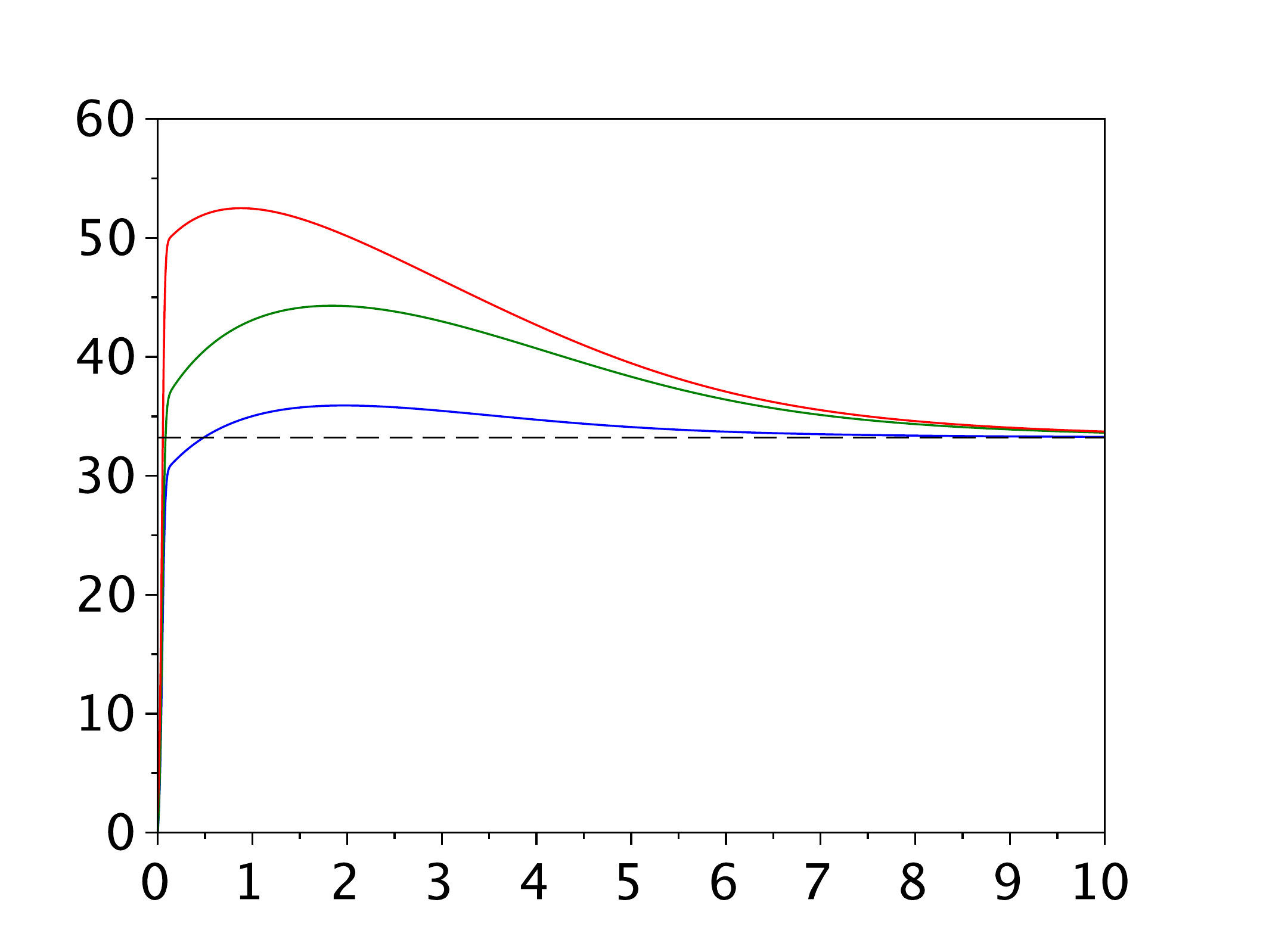}
\includegraphics[width=6cm]{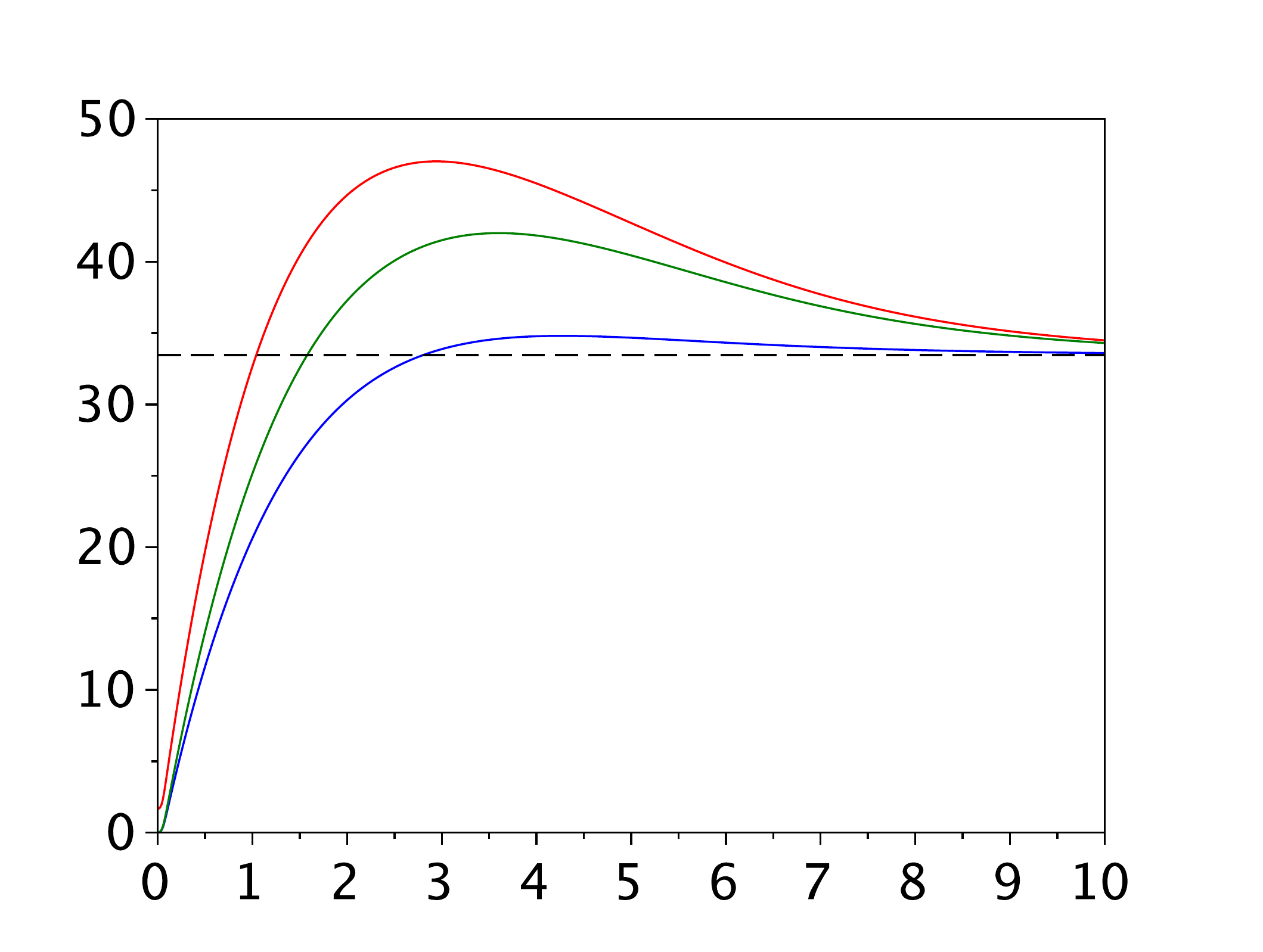}
\caption{{\bf Release of larvae}: evolution of the uninfected (top) and {\em Wolbachia}-infected (bottom), as function of time.
The larvae appear on the left column, the adults on the right one.
The components of the state $x$ (resp.\ of the estimates $x_-$, $x^+$) appear in {\color{green}\bf green} (resp.\ in {\color{blue}\bf blue}, in {\color{red}\bf red}).}
\label{fig2}
\end{figure*}

\begin{theorem}\em (Sufficient conditions for successful introduction via release of adults).
\label{th4}
Assume available upper and lower state estimates $x_\pm$ for the solution of \eqref{eq1}, according to \eqref{eq19} and \eqref{eq28}.
Assume that for some $T\geq 0$, \eqref{eq8b} holds for given constants $k>0$, $k_U>\cR_0^U-1$.
Then
\vspace{-.1cm}
\addtocounter{equation}{-1}
\begin{equation}
\label{eq9}
\hspace{-.3cm}
\lim_{t\to +\infty} \begin{pmatrix}
L_U(t) \\ A_U(t)
\end{pmatrix}
= \begin{pmatrix}
0 \\ 0
\end{pmatrix},\
\liminf_{t\to +\infty} \begin{pmatrix}
L_W(t) \\ A_W(t)
\end{pmatrix} \geq \begin{pmatrix}
L_W^* \\ A_W^*
\end{pmatrix}
\end{equation}
\end{theorem}

\vspace{.2cm}
Formula \eqref{eq8b} is displayed on top of the next page.
The goal of the control is to steer the system to the complete infestation equilibrium $(0,0,L_W^*,A_W^*)$.
Notice that when it succeeds, the lower bound of $u_A$ in \eqref{eq8b} may vanish asymptotically {\em only if}
$k = \gamma_W$.
In this case, the second inequality in \eqref{eq8b} becomes (in view of \eqref{eq12})
\vspace{-.2cm}
\begin{equation}
\label{eq32}
\hspace{-.3cm}
u_A \geq 
\begin{pmatrix}
k_U & k_U(\gamma_W -\gamma_U) & -1
\end{pmatrix}
\begin{pmatrix}
L_U^+\\
A_U^+\\
L_W^-
\end{pmatrix}
\text{ on } [T,+\infty)
\end{equation}
expression in which the three useful components of the estimator pertain to $x_-$, see \eqref{eq28}.

{\em Proof of Theorem {\rm\ref{th4}}.}

\noindent $\bullet$ 1.\
We have
$\dot A_W-k_U\dot A_U
= L_W -\gamma_W  A_W + u_A-k_U(L_U -\gamma_U  A_U)$.
By assumption, one has on $[T,+\infty)$:
$u_A \geq - L_W^- +k_UL_U^+ + |\gamma_W -k|_+ A_W^+ - |\gamma_W -k|_- A_W^-
 +k_U|k -\gamma_U|_+A_U^+
-k_U|k -\gamma_U|_-A_U^-$.
Thus
\vspace{-.1cm}
\begin{eqnarray*}
\lefteqn{\dot A_W-k_U\dot A_U}\\
& \geq &
L_W - L_W^-+k_U(L_U^+-L_U)\\
& &
-\gamma_W  A_W
 + |\gamma_W -k|_+ A_W^+
 - |\gamma_W -k|_- A_W^-\\
 & &
+k_U\left(
\gamma_U  A_U
+|k -\gamma_U|_+A_U^+
-|k -\gamma_U|_-A_U^-
\right)\\
& \geq &
-\gamma_W  A_W
 + |\gamma_W -k|_+ A_W^+
 - |\gamma_W -k|_- A_W^-\\
 & &
 +k_U\left(
 \gamma_U  A_U
+|k -\gamma_U|_+A_U^+
-|k -\gamma_U|_-A_U^-
\right)
\end{eqnarray*}

Now, using \eqref{eq26}, one verifies that
$-\gamma_W  A_W
 + |\gamma_W -k|_+ A_W^+
 - |\gamma_W -k|_- A_W^-
 \geq
 -kA_W$,
 and similarly
$\gamma_U  A_U
+|k -\gamma_U|_+A_U^+
-|k -\gamma_U|_-A_U^-
\geq
kA_U$.
One thus deduces from these two inequalities the key property:
\begin{equation*}
\dot A_W-k_U \dot A_U \geq -k(A_W-k_UA_U)\ .
\end{equation*}
Integrating this differential inequality yields, for any $t\geq T$,
$A_W(t)-k_UA_U(t) \geq e^{-k(t-T)} (A_W(T)-k_UA_U(T))$.
Thus
\begin{equation}
\label{eq20}
\liminf_{t\to +\infty} (A_W(t)-k_UA_U(t)) \geq 0\ .
\end{equation}

\noindent $\bullet$ 2.\
As $k_U>\cR_0^U-1$, there exists $\varepsilon >0$ such that $k_U>\cR_0^U(1+\varepsilon)-1$.
Due to \eqref{eq20}, there exists $T'\geq T$ such that, for any $t\geq T'$,
$A_W(t)\geq (\cR_0^U(1+\varepsilon)-1)A_U(t)$,
and, due to the fact that the function $\Rset_+\to\Rset_+$, $z\mapsto \frac{1}{1+z}$ is decreasing,
$\frac{A_U(t)}{A_U(t)+A_W(t)} \leq \frac{1}{1+\cR_0^U(1+\varepsilon)-1} = \frac{1}{\cR_0^U(1+\varepsilon)}$.
For $t\geq T'$, the evolution of the unifected therefore obeys the following inequalities:
$\dot L_U \leq \frac{\gamma_U}{1+\varepsilon} A_U - (1+L_W+L_U)L_U
\leq \frac{\gamma_U}{1+\varepsilon} A_U - (1+L_U)L_U$,
$\dot A_U = L_U-\gamma_UA_U$.
As system \eqref{eq13} with $\gamma:=\gamma_U$, $\cR := \frac{1}{1+\varepsilon}<1$ is monotone, it may serve as a comparison system for the subsystem describing the evolution of $(L_U,A_U)$: denoting $(L,A)$ the solution of system \eqref{eq13} such that $L(T')=L_U(T')$, $A(T')=A_U(T')$, one has, for any $t\geq T'$,
$0 \leq L_U(t) \leq L(t)$, $0 \leq A_U(t) \leq A(t)$.
On the other hand, as $\cR<1$, $(L(t),A(t))$ vanishes when $t\to +\infty$.
One thus gets:
\vspace{-.1cm}
\begin{equation}
\label{eq22}
\lim_{t\to +\infty} L_U(t) = \lim_{t\to +\infty} A_U(t) = 0\ ,
\end{equation}
which is the first inequality of \eqref{eq9}.


\noindent $\bullet$ 3.\
We now deduce lower bounds on $L_W,A_W$ from the identities in \eqref{eq22}.
For any $\varepsilon>0$, there exists $T''>T$ such that, for any $t\geq T''$, $0\leq L_U(t)\leq\varepsilon$.
Using the fact that $u_A,u_L$ have nonnegative values, one gets:
$\dot L_W \geq \gamma_W\cR_0^W A_W - (1+\varepsilon+L_W)L_W$,
$\dot A_W \geq L_W-\gamma_WA_W$.
The underlying system
\begin{equation}
\label{eq24}
\dot L' = \gamma_W\cR_0^W A' - (1+\varepsilon+L')L',\quad
\dot A' = L'-\gamma_WA'
\end{equation}
is monotone.
Applying again comparison result \cite{Coppel:1965aa}, now to the subsystem $(L_W,A_W)$, one gets that, for any $t\geq T''$,
\begin{equation}
\label{eq25}
L_W(t) \geq L'(t),\quad A_W(t) \geq A'(t)
\end{equation}
for the trajectory of \eqref{eq24} initialized by $L'(T'')=L_W(T'')$, $A'(T'')=A_W(T'')$.

Now, \eqref{eq24} is a variant of \eqref{eq13}, and may be analyzed in the same way.
It is not difficult to show that when $1+\varepsilon<\cR_0^W$, every nonzero trajectory of \eqref{eq24} converges towards the positive equilibrium $(\cR_0^W-1-\varepsilon,\frac{\cR_0^W-1-\varepsilon}{\gamma_W})$.
For $\varepsilon>0$ such that $\varepsilon<\cR_0^W-1$, one thus has
${\displaystyle\lim_{t\to +\infty}} L'(t) = \cR_0^W-1-\varepsilon$, ${\displaystyle\lim_{t\to +\infty}} A'(t) = \frac{\cR_0^W-1-\varepsilon}{\gamma_W}$.
Deducing from \eqref{eq25} that
${\displaystyle\liminf_{t\to +\infty}} L_W(t) \geq {\displaystyle\lim_{t\to +\infty}} L'(t) = \cR_0^W-1-\varepsilon$, 
${\displaystyle\liminf_{t\to +\infty}} A_W(t) \geq {\displaystyle\lim_{t\to +\infty}} A'(t) = \frac{\cR_0^W-1-\varepsilon}{\gamma_W}$
and passing to the limit $\varepsilon\to 0$, yields the second inequality in \eqref{eq9}.
This completes the proof of Theorem \ref{th4}.
\hfill $\blacksquare$

\subsection{Control by release of larvae}
\label{se42}

We now examine a control approach based on release of larvae.
This method amounts to acting on the input term $u_L$ in \eqref{eq1c}.
The only hindrance to the growth of the infected population lies on the competition term in the previous equation.
Therefore, if the rate of introduction of infected larvae is sufficient to compensate for the competition effect, then the infected population should reach the complete infestation equilibrium and induce collapse of the uninfected population, through the competitive pressure term in \eqref{eq1b}.
This is evidenced by the following result.


\begin{theorem}\em (Sufficient conditions for successful introduction via release of larvae).
\label{th5}
Assume available upper and lower state estimates $x_\pm$ for the solution of \eqref{eq1}, according to \eqref{eq19} and \eqref{eq28}.
Assume that for some $T\geq 0$,
\vspace{-.1cm}
{\renewcommand{\theequation}{\bf L}
\begin{equation}
\label{eq8a}
u_L > L_U^+L_W^+,\quad
u_A \geq 0 \qquad
\text{ on } [T,+\infty)\ .
\end{equation}
}
Then \eqref{eq9} holds.
\end{theorem}

\vspace{.1cm}
{\em Proof of Theorem {\rm\ref{th5}}.}

\noindent $\bullet$ 1.\
When \eqref{eq8a} is fulfilled, then the infected population obeys the following differential inequalities:
\addtocounter{equation}{-1}
$\dot L_W > \gamma_W {\cal R}_0^W  A_W - (1+ L_W+ L_U) L_W + L_U^+L_W^+
\geq \gamma_W {\cal R}_0^W  A_W - (1+ L_W) L_W$
and
$\dot A_W \geq L_W -\gamma_W  A_W$.
One may thus bound from below the solutions of the subsystem $(L_W,A_W)$ by the solution $(L,A)$ of system \eqref{eq13} with $\gamma:=\gamma_W$ and $\cR :=\cR_0^W$, initialised with the same valor than $(L_W,A_W)$ at a certain $T'\geq T$.
This solution of \eqref{eq13} tends to the positive equilibrium, or to the zero equilibrium {\em only if the initial condition is precisely the zero equilibrium.}

\noindent $\bullet$ 2.\
We now use the first inequality in \eqref{eq8a}.
Due to its {\em strict} character\footnote{ The use of $>$ in this context is defined in the Notations, see Section \ref{se1}.}, any trajectory of \eqref{eq1} takes on some {\em nonzero} value on $[T,+\infty)$, say at time $T'$.
As $(0,0)< (L(T'),A(T'))$, convergence of $(L,A)$ to the zero equilibrium of \eqref{eq13} is excluded, and this solution is forced  to converge towards the positive equilibrium.
Using this property one deduces that the second inequality in \eqref{eq9} is fulfilled.

\noindent $\bullet$ 3.\
In turn, we deduce from the preceding inference, that the evolution of the non-infected population fulfills asymptotically the inequalities:
\vspace{-.1cm}
\begin{multline}
\label{eq11}
\dot L_U \leq \gamma_U{\cal R}_0^U \frac{ A_U}{ A_U+ A_W^*} A_U - (1+ L_W^*+ L_U) L_U,\\
\dot A_U \leq  L_U -\gamma_U  A_U
\end{multline}
\vspace{-.1cm}
The underlying differential system
\begin{equation}
\label{eq29}
\hspace{-.1cm}
\dot L = \gamma_U{\cal R}_0^U \frac{ A}{ A+ A_W^*} A - (1+ L_W^*+ L) L,\
\dot A =  L -\gamma_U  A
\end{equation}
is monotone, admits only bounded trajectories and possesses $(0,0)$ as unique equilibrium (otherwise, a supplementary equilibrium would exist for the initial system \eqref{eq1}).
Therefore \cite[Theorem 3.1, p.\ 18]{Smith:1995aa}, every trajectory of system \eqref{eq29} converges towards the unique equilibrium $(0,0)$.
Applying comparison result to subsystem $(L_U,A_U)$ shows that the system of differential inequalities \eqref{eq11} has nonnegative solutions that are at most equal to those of \eqref{eq29}.
One infers the first inequality of \eqref{eq9}, and this achieves the proof of Theorem \ref{th5}.
\hfill $\blacksquare$

\section{Numerical simulations}
\label{se6}

Numerical simulations of the controlled system \eqref{eq00}-\eqref{eq18} have been achieved for the control laws designed in Theorems \ref{th4} and \ref{th5} (using the {\tt stiff} option of the {\tt ode} function of the free open-source software Scilab).
The following numerical values, extracted from \cite{Bliman:2018aa}, have been adopted:
$\gamma_U = 0.79365$, $\gamma_W = 0.99207$, $\cR_0^U = 45$, $\cR_0^W = 34.2$, so that
$L_U^*= 44.0$, $A_U^*= 55.4$, $L_W^*= 33.2$, $A_W^*=33.5$.

The initial state is systematically taken at the disease-free equilibrium $(L_U^*,A_U^*,0,0)$.
For all simulations, observer initialization is done with quite conservative estimates:
$(L_U^-(0),A_U^-(0),L_W^-(0),A_W^-(0)) = (0,0,0,0)$
and $(L_U^+(0),A_U^+(0),L_W^+(0),A_W^+(0)) = (2L_U^*,2A_U^*,0.05L_W^*,0.05A_W^*)$.
For each type of release, noisy measurement is considered, with $y_-=80\% y$, $y_+=120\% y$.

The choice of the gains present no specific difficulty, as the conditions on the observer gains stated in Theorem \ref{th0} are not hard to meet.
The output gain $C$ is taken as the example in Remark \ref{re0}, while the observer gains $K_\pm(\cdot)$ are chosen as explained in Remark \ref{re1}, with
\vspace{-.1cm}
\[
\varepsilon = 10^{-5},\qquad
M_- = M_+ = 0.1
\begin{pmatrix}
-1 & -1 & 1 & 1\\
-1 & -1 & 1 & 1
\end{pmatrix}^{\mbox{\tiny\sf T}}\ .
\]
In these conditions, all assumptions of Theorem \ref{th2} are fulfilled.
The simulations show, as expected and in spite of imprecise initial state estimation and measurement, the convergence of the controlled system towards the complete infestation equilibrium $(0,0,L_W^*,A_W^*)$, in dashed line in the two Figures.

\paragraph{Release of adult mosquitos (Theorem \ref{th4})}

The control $u:=(u_L,u_A)$ is defined by taking equalities in the formulas in \eqref{eq8b}.
One chooses  $k_U=1.1(\cR_0^U-1)$ and $k=\gamma_W$.
The control $u_A$ is then linear\footnote{ More precisely:
$u_A = \left|
\begin{pmatrix}
k_U & k_U(\gamma_W -\gamma_U) & -1
\end{pmatrix}
\begin{pmatrix}
L_U^+\\
A_U^+\\
L_W^-
\end{pmatrix}
\right|_+$.}, equal to the right-hand side of \eqref{eq32}.
As noticed after the statement of Theorem \ref{th4}, only $x_-$ is used for the control synthesis.
The results are shown in Figure \ref{fig1}.
The green curves represent the `true' values (components of $x$), while the blue ones show the evolution of the components $L_U^+, A_U^+, L_W^-, A_W^-$ of $x_-$.

\paragraph{Release of larvae (Theorem \ref{th5})}

The control $u:=(u_L,u_A)$ is defined here by taking equalities in \eqref{eq8a}.
Notice that the strict inequality in \eqref{eq8a} is obtained thanks to a nonzero initialization of $L_W^+$.
The results are shown in Figure \ref{fig2}.
The green curves represent the components of $x$, while the blue (resp.\ red) ones show the components of $x_-$ (resp.\ $x^+$).


\section{Concluding remarks}
\label{se7}

Two release strategies have been proposed, allowing to realize complete infestation of a population of {\em Aedes} mosquitoes by a population artificially infected by the bacterium {\em Wolbachia}, which reduces their ability as vectors of several important diseases.
The efficiency of the approach 
has been demonstrated analytically and tested numerically.
Its main force is to be expressed in general terms, giving rise to potential extensions to other models in the literature, e.g.\ those in \cite{Keeling:2003aa,Farkas:2010aa,Hughes:2013aa,Zheng:2014aa,Koiller:2014aa,Ndii:2015aa,Zheng:2017aa,Yakob:2017aa,Xue:2017aa,Almeida:2018aa}, due to their monotonicity\footnote{ One may check that the underlying entomological models are all monotone, assuming where appropriate complete cytoplasmic incompatibility, zero delay and possibly after merging of some male/female compartments.}.
Being based on monotone system properties, it is believed to be a good way to ensure robust behavior with respect to unmodeled dynamics and parameter uncertainties, with no need of precise knowledge of the boundary between the basins of attraction of the two locally asymptotically stable equilibria.
Future works include the consideration of measurement delays.
Also, further study is needed in order to reduce the peaking effect apparent with the first method.
%
%
Last, extensions to non monotone systems should be considered, as well as more realistic impulsive periodic releases, much in the spirit of \cite{Bliman:2018ab} in the context of Sterile Insect Techniques.

\bibliographystyle{IEEEtran}
\bibliography{IEEEabrv,ecc2019}

\begin{thebibliography}{10}
\providecommand{\url}[1]{#1}
\csname url@samestyle\endcsname
\providecommand{\newblock}{\relax}
\providecommand{\bibinfo}[2]{#2}
\providecommand{\BIBentrySTDinterwordspacing}{\spaceskip=0pt\relax}
\providecommand{\BIBentryALTinterwordstretchfactor}{4}
\providecommand{\BIBentryALTinterwordspacing}{\spaceskip=\fontdimen2\font plus
\BIBentryALTinterwordstretchfactor\fontdimen3\font minus
  \fontdimen4\font\relax}
\providecommand{\BIBforeignlanguage}[2]{{%
\expandafter\ifx\csname l@#1\endcsname\relax
\typeout{** WARNING: IEEEtran.bst: No hyphenation pattern has been}%
\typeout{** loaded for the language `#1'. Using the pattern for}%
\typeout{** the default language instead.}%
\else
\language=\csname l@#1\endcsname
\fi
#2}}
\providecommand{\BIBdecl}{\relax}
\BIBdecl

\bibitem{Bliman:2018aa}
P.-A. Bliman, M.~S. Aronna, F.~C. Coelho, and M.~A. Da~Silva, ``Ensuring
  successful introduction of {Wolbachia in natural populations of Aedes
  aegypti} by means of feedback control,'' \emph{Journal of mathematical
  biology}, vol.~76, no.~5, pp. 1269--1300, 2018.

\bibitem{Alphey:2010aa}
L.~Alphey, M.~Benedict, R.~Bellini, G.~G. Clark, D.~A. Dame, M.~W. Service, and
  S.~L. Dobson, ``{Sterile-Insect Methods for Control of Mosquito-Borne
  Diseases: An Analysis},'' \emph{Vector-Borne and Zoonotic Diseases}, vol.~10,
  no.~3, pp. 295--311, apr 2010.

\bibitem{Alphey:2014aa}
L.~Alphey, ``Genetic control of mosquitoes,'' \emph{Annual review of
  entomology}, vol.~59, 2014.

\bibitem{McMeniman:2009aa}
C.~J. McMeniman, R.~V. Lane, B.~N. Cass, A.~W. Fong, M.~Sidhu, Y.-F. Wang, and
  S.~L. O'neill, ``Stable introduction of a life-shortening {Wolbachia}
  infection into the mosquito {Aedes} aegypti,'' \emph{Science}, vol. 323, no.
  5910, pp. 141--144, 2009.

\bibitem{Moreira:2009aa}
L.~A. Moreira, I.~Iturbe-Ormaetxe, J.~A. Jeffery, G.~Lu, A.~T. Pyke, L.~M.
  Hedges, B.~C. Rocha, S.~Hall-Mendelin, A.~Day, M.~Riegler, L.~E. Hugo, K.~N.
  Johnson, B.~H. Kay, E.~A. McGraw, A.~F. van~den Hurk, P.~A. Ryan, and S.~L.
  O'Neill, ``{A Wolbachia Symbiont in Aedes aegypti Limits Infection with
  Dengue, Chikungunya, and Plasmodium},'' \emph{Cell}, vol. 139, no.~7, pp.
  1268--1278, dec 2009.

\bibitem{Hoffmann:2011aa}
A.~A. Hoffmann, B.~L. Montgomery, J.~Popovici, I.~Iturbe-Ormaetxe, P.~H.
  Johnson, F.~Muzzi, M.~Greenfield, M.~Durkan, Y.~S. Leong, Y.~Dong, H.~Cook,
  J.~Axford, A.~G. Callahan, N.~Kenny, C.~Omodei, E.~A. McGraw, P.~A. Ryan,
  S.~A. Ritchie, M.~Turelli, and S.~L. O'Neill, ``{Successful establishment of
  Wolbachia in Aedes populations to suppress dengue transmission},''
  \emph{Nature}, vol. 476, no. 7361, pp. 454--457, aug 2011.

\bibitem{Wilder-Smith:2017aa}
A.~Wilder-Smith, D.~J. Gubler, S.~C. Weaver, T.~P. Monath, D.~L. Heymann, and
  T.~W. Scott, ``Epidemic arboviral diseases: priorities for research and
  public health,'' \emph{The Lancet infectious diseases}, vol.~17, no.~3, pp.
  e101--e106, 2017.

\bibitem{Keeling:2003aa}
M.~J. Keeling, F.~Jiggins, and J.~M. Read, ``The invasion and coexistence of
  competing {Wolbachia strains},'' \emph{Heredity}, vol.~91, no.~4, p. 382,
  2003.

\bibitem{Farkas:2010aa}
J.~Farkas and P.~Hinow, ``Structured and unstructured continuous models for
  {Wolbachia} infections,'' \emph{Bulletin of mathematical biology}, vol.~72,
  no.~8, pp. 2067--2088, 2010.

\bibitem{Zheng:2014aa}
B.~Zheng, M.~Tang, and J.~Yu, ``Modeling {Wolbachia} spread in mosquitoes
  through delay differential equations,'' \emph{SIAM Journal on Applied
  Mathematics}, vol.~74, no.~3, pp. 743--770, 2014.

\bibitem{Zheng:2017aa}
B.~Zheng and Y.~Xiao, ``A mathematical model to assess the effect of the
  constant release policy on population suppression,'' \emph{Nonlinear Analysis
  and Differential Equations}, vol.~5, no.~4, pp. 197--207, 2017.

\bibitem{Yakob:2017aa}
L.~Yakob, S.~Funk, A.~Camacho, O.~Brady, and W.~J. Edmunds, ``Aedes aegypti
  control through modernized, integrated vector management,'' \emph{PLoS
  currents}, vol.~9, 2017.

\bibitem{Xue:2017aa}
L.~Xue, C.~A. Manore, P.~Thongsripong, and J.~M. Hyman, ``Two-sex mosquito
  model for the persistence of {Wolbachia},'' \emph{Journal of biological
  dynamics}, vol.~11, no. sup1, pp. 216--237, 2017.

\bibitem{Campo-Duarte:2017aa}
D.~E. Campo-Duarte, D.~Cardona-Salgado, and O.~Vasilieva, ``Establishing
  {wMelPop Wolbachia} infection among wild {Aedes} aegypti females by optimal
  control approach,'' \emph{Appl Math Inf Sci}, vol.~11, no.~4, pp. 1011--1027,
  2017.

\bibitem{Campo-Duarte:2018aa}
D.~E. Campo-Duarte, O.~Vasilieva, D.~Cardona-Salgado, and M.~Svinin, ``Optimal
  control approach for establishing {wMelPop Wolbachia} infection among wild
  {Aedes} aegypti populations,'' \emph{Journal of mathematical biology},
  vol.~76, no.~7, pp. 1907--1950, 2018.

\bibitem{Almeida:2018aa}
\BIBentryALTinterwordspacing
L.~Almeida, Y.~Privat, M.~Strugarek, and N.~Vauchelet, ``{Optimal releases for
  population replacement strategies, application to Wolbachia},'' Jun. 2018,
  working paper or preprint. [Online]. Available:
  \url{https://hal.archives-ouvertes.fr/hal-01807624}
\BIBentrySTDinterwordspacing

\bibitem{Hughes:2013aa}
H.~Hughes and N.~F. Britton, ``{Modelling the use of Wolbachia to control
  dengue fever transmission},'' \emph{Bulletin of mathematical biology},
  vol.~75, no.~5, pp. 796--818, 2013.

\bibitem{Koiller:2014aa}
J.~Koiller, M.~Da~Silva, M.~Souza, C.~Code{\c c}o, A.~Iggidr, and G.~Sallet,
  ``{Aedes, Wolbachia and Dengue},'' Inria, France, Research Report RR-8462,
  Jan. 2014.

\bibitem{Ndii:2015aa}
M.~Z. Ndii, R.~I. Hickson, D.~Allingham, and G.~Mercer, ``Modelling the
  transmission dynamics of dengue in the presence of {Wolbachia},''
  \emph{Mathematical biosciences}, vol. 262, pp. 157--166, 2015.

\bibitem{Focks:2004aa}
D.~A. Focks \emph{et~al.}, ``A review of entomological sampling methods and
  indicators for dengue vectors,'' Geneva: World Health Organization, Tech.
  Rep., 2004.

\bibitem{Silver:2007aa}
J.~B. Silver, \emph{Mosquito ecology: field sampling methods}.\hskip 1em plus
  0.5em minus 0.4em\relax Springer Science \& Business Media, 2007.

\bibitem{Almeida:2019aa}
L.~Almeida, M.~Duprez, Y.~Privat, and N.~Vauchelet, ``Control strategies on
  mosquitos population for the fight against arboviruses,'' \emph{arXiv
  preprint arXiv:1901.05688}, 2019.

\bibitem{Moisan:2009aa}
M.~Moisan, O.~Bernard, and J.-L. Gouz{\'e}, ``Near optimal interval observers
  bundle for uncertain bioreactors,'' \emph{Automatica}, vol.~45, no.~1, pp.
  291--295, 2009.

\bibitem{Efimov:2013aa}
D.~Efimov, T.~Ra{\"\i}ssi, S.~Chebotarev, and A.~Zolghadri, ``Interval state
  observer for nonlinear time varying systems,'' \emph{Automatica}, vol.~49,
  no.~1, pp. 200--205, 2013.

\bibitem{Efimov:2016aa}
D.~Efimov and T.~Ra{\"\i}ssi, ``Design of interval observers for uncertain
  dynamical systems,'' \emph{Automation and Remote Control}, vol.~77, no.~2,
  pp. 191--225, 2016.

\bibitem{Yang:2009aa}
H.~Yang, M.~d. L. d.~G. Macoris, K.~Galvani, M.~Andrighetti, and D.~Wanderley,
  ``Assessing the effects of temperature on the population of {Aedes} aegypti,
  the vector of dengue,'' \emph{Epidemiology \& Infection}, vol. 137, no.~8,
  pp. 1188--1202, 2009.

\bibitem{Ferreira:2014aa}
C.~P. Ferreira and W.~A. Godoy, \emph{Ecological modelling applied to
  entomology}.\hskip 1em plus 0.5em minus 0.4em\relax Springer, 2014.

\bibitem{Yeap:2011aa}
H.~L. Yeap, P.~Mee, T.~Walker, A.~R. Weeks, S.~L. O'Neill, P.~Johnson, S.~A.
  Ritchie, K.~M. Richardson, C.~Doig, N.~M. Endersby \emph{et~al.}, ``Dynamics
  of the ``popcorn'' {Wolbachia} infection in outbred {Aedes} aegypti informs
  prospects for mosquito vector control,'' \emph{Genetics}, vol. 187, no.~2,
  pp. 583--595, 2011.

\bibitem{Ruang-areerate:2006aa}
T.~Ruang-Areerate and P.~Kittayapong, ``Wolbachia transinfection in {Aedes}
  aegypti: a potential gene driver of dengue vectors,'' \emph{Proceedings of
  the National Academy of Sciences}, vol. 103, no.~33, pp. 12\,534--12\,539,
  2006.

\bibitem{Southwood:1972aa}
T.~Southwood, G.~Murdie, M.~Yasuno, R.~J. Tonn, and P.~Reader, ``Studies on the
  life budget of {Aedes aegypti in Wat Samphaya, Bangkok, Thailand},''
  \emph{Bulletin of the World Health Organization}, vol.~46, no.~2, p. 211,
  1972.

\bibitem{Focks:1993aa}
D.~A. Focks, D.~Haile, E.~Daniels, and G.~A. Mount, ``Dynamic life table model
  for {Aedes aegypti (Diptera: Culicidae): analysis of the literature and model
  development},'' \emph{Journal of medical entomology}, vol.~30, no.~6, pp.
  1003--1017, 1993.

\bibitem{Influential:1998aa}
S.~L. {O'Neill}, A.~A. Hoffman, and J.~H. Werren, Eds., \emph{{Influential
  Passengers: Inherited Microorganisms and Arthropod Reproduction}}.\hskip 1em
  plus 0.5em minus 0.4em\relax Oxford University Press, 1998.

\bibitem{Walker:2011aa}
T.~Walker, P.~Johnson, L.~Moreira, I.~Iturbe-Ormaetxe, F.~Frentiu,
  C.~McMeniman, Y.~S. Leong, Y.~Dong, J.~Axford, P.~Kriesner \emph{et~al.},
  ``{The {\em w}Mel Wolbachia strain blocks dengue and invades caged Aedes
  aegypti populations},'' \emph{Nature}, vol. 476, no. 7361, pp. 450--453,
  2011.

\bibitem{Smith:1995aa}
H.~L. Smith, \emph{Monotone Dynamical Systems. An Introduction to the Theory of
  Competitive and Cooperative Systems}, ser. Mathematical Surveys and
  Monographs.\hskip 1em plus 0.5em minus 0.4em\relax American Mathematical
  Society, 1995, vol.~41.

\bibitem{Coppel:1965aa}
W.~A. Coppel, \emph{Stability and asymptotic behavior of differential
  equations}, ser. Heath mathematical monographs.\hskip 1em plus 0.5em minus
  0.4em\relax Heath, 1965.

\bibitem{Bliman:2018ab}
P.-A. Bliman, D.~Cardona-Salgado, Y.~Dumont, and O.~Vasilieva, ``Implementation
  of control strategies for sterile insect techniques,'' \emph{arXiv preprint
  arXiv:1812.01277}, 2018.

\end{thebibliography}


\end{document}